\pdfoutput=1
\RequirePackage{ifpdf}
\ifpdf 
\documentclass[pdftex]{sigma}
\else
\documentclass{sigma}
\fi

\newcommand{\R}{\mathbb{R}}
\newcommand{\C}{\mathbb{C}}

\DeclareMathOperator{\sgn}{sgn}
\DeclareMathOperator{\spn}{span}

\begin{document}

\allowdisplaybreaks

\renewcommand{\PaperNumber}{011}

\FirstPageHeading

\ShortArticleName{Porous Medium Equation and~Global Actions of the Symmetry Group}

\ArticleName{On the $\boldsymbol{n}$-Dimensional Porous Medium Dif\/fusion
\\
Equation and Global Actions of the Symmetry Group}

\Author{Jose A.~FRANCO}

\AuthorNameForHeading{J.A.~Franco}

\Address{Department of Mathematics and~Statistics, University of North Florida,
\\
1 UNF Drive, Jacksonville, FL 32224 USA}
\Email{\href{mailto:jose.franco@unf.edu}{jose.franco@unf.edu}}

\ArticleDates{Received September 10, 2012, in f\/inal form February 08, 2013; Published online February 12, 2013}

\Abstract{By restricting to a~special class of smooth functions, the local action of the symmetry
group is globalized.
This special class of functions is constructed using parabolic induction.}


\Classification{22E70; 35Q35}

\section{Introduction}
The theory of Lie groups f\/inds its genesis when Sophus Lie sets up the task to develop an
analogue of Galois theory for dif\/ferential equations.
The original prolongation algorithm of Sophus Lie provides a~set of point symmetries of
a~dif\/ferential equation (for a~modern treatment, see~\cite{MR1240056}).
It is well-known that these point symmetries generate a~\textit{local} Lie group action on the
space of solutions of the dif\/ferential equation and~this action seldom globalizes.
However, most of the modern results on Lie group theory apply to \textit{global} Lie groups
and~\textit{global} Lie group representations.
As a~result, many of the standard techniques of representation theory are not always applicable to
dif\/ferential equations.
Therefore, it is an important problem to f\/ind a~globalization of the local action of the symmetry
group of a~dif\/ferential equation (see~\cite{Craddock2, Franco, Sepanski3, Kostant, Sepanski}
and~references therein).
This is usually achieved by
restricting to a~special class of functions.

In this paper, we study the globalization problem for the $n$-dimensional porous medium equation
\begin{gather}
\label{Porous}
u_t=\Delta_n(u^m),
\end{gather}
where $\Delta_n$ is the $n$-dimensional Laplacian and~$u$ is a~function of
$x=(x_1,\dots,x_n)\in\R^n$ and~$t\in\R$.
When $m>1$, this equation models slow dif\/fusion phenomena so this condition is often assumed
(cf.~\cite{wu2001nonlinear}).
However, for reasons that will become evident below, we will allow $m\in\R\setminus\{0,1\}$.
It is worth noting that the globalization problem for the special case $m=1$ was solved
in~\cite{Sepanski}, where applications of the global action of the group can be found.
We will pay particular attention to the case $m=\frac{n-2}{n+2}$, because for this special value of
$m$ the symmetry group is $n$-dimensions larger than in the generic case.

The goal of this paper is to describe a~class of functions on which the action of the symmetry
group globalizes and~to describe the action of the group on these functions.
This idea is quite prof\/itable when applied to linear PDE's (see~\cite{Franco, Sepanski3,
Sepanski}).
However, it can also be applied to nonlinear equations as well.
For example, the globalization problem was studied for a~family of nonlinear heat equations
in~\cite{Sepanski2} and~for the nonlinear potential f\/iltration equation in~\cite{Sepanski4}.
In both cases, the group action is globalized by using parabolic induction on solvable groups.
In this article, the globalization is achieved by parabolic induction on a~semisimple group.
In~\cite{Sepanski2, Sepanski4}, the global actions of the group are not linear.
In contrast, the globalization of the action of the symmetry group of the porous medium equation is
linear.

The paper is organized as follows.
In Section~\ref{S.G.}, we list the inf\/initesimal generators of the symmetry group of
equation~\eqref{Porous}.
The symmetry group for the
one-dimensional case is calculated in~\cite{Ovsiannikov} (see also~\cite{CRC,dos2011enhanced}
and~references therein).
For the $n$-dimensional case, it was be found in~\cite{Dorodnitsyn} (the
result is available in the handbook~\cite{CRC}).
In either case, the symmetry group is realized as a~subgroup of
$G:=\text{SL}(2,\R)\times\text{SO}(n+1,1)_0$, where $\text{SO}(n+1,1)_0$ stands for the identity
component of the generalized special orthogonal group with signature $(n+1,1)$.

In Section~\ref{Three}, we brief\/ly study the structure of the group $G$.
In Section~\ref{induced}, we use parabolic induction to construct a~family of representations
of~$G$ and~by restriction, of the symmetry group of~\eqref{Porous}.
We use these representations to def\/ine class a~of smooth functions on which the action of the
symmetry group globalizes (see equation~\eqref{Iprime}).
In Section~\ref{actions}, we write explicitly the action of the $1$-parameter subgroups generated
by the inf\/initesimal generators of~$G$.
Finally, in Section~\ref{applications}, we realize well-known solutions of equation~\eqref{Porous} as
elements of the constructed representation of~$G$ and~we point references for possible
generalizations of this work.

\section{Symmetry group}
\label{S.G.}

Using Lie's prolongation algorithm, we calculate the symmetry group of equation~\eqref{Porous} and~re-obtain the result of~\cite{Dorodnitsyn}.
This yields two cases depending on the value of $m$.
When $m
\neq\frac{n-2}{n+2}$, the inf\/initesimal generators of the symmetry group are
\begin{gather}
X_1=\partial_t,\label{GenGenerators}
\qquad
X_2=\sum_{i=1}^nx_i\partial_i+\frac{2u}{m-1}\partial_u,
\qquad
X_3=-t\partial_t+\frac{u}{m-1}\partial_u,
\\
Y_i=\partial_i\qquad\text{for}\quad 1\leq i\leq n,
\\[1.8mm]
Z_{i,j}=x_i\partial_j-x_j\partial_i\qquad\text{for}\quad 1\leq i<j\leq n.\label{LastGen}
\end{gather}
Let $\mathfrak{s}$ denote the parabolic subalgebra of upper triangular matrices in
$\mathfrak{sl}_2(\R)$ and~let $\mathfrak{g}:=\mathfrak{s}\oplus\mathfrak{so}(n+1,1)$.
Then, the inf\/initesimal generators~\eqref{GenGenerators}--\eqref{LastGen} span an algebra
isomorphic to a~parabolic subalgebra of $\mathfrak g$.

In addition to these generators, in the special case $m=\frac{n-2}{n+2}$ the symmetry group is
extended by the one-parameter groups generated by the operators
\begin{gather}
\label{ExtraGen}
W_i=\left(x_i^2-\sum_{j\neq i}x_j^2\right)\partial_i+\sum_{j\neq i}2x_ix_j\partial_j+\frac{4x_i u}{m-1}\partial_u
\end{gather}
for $1\leq i\leq n$.
The inf\/initesimal generators~\eqref{GenGenerators}--\eqref{ExtraGen} span an algebra isomorphic
to
$\mathfrak{g}$.
The isomorphism for the $\mathfrak{so}(n+1,1)$ part is explicitly def\/ined in the following way
\begin{gather*}
X_2\mapsto-E_{n+1,n+2}-E_{n+2,n+1},
\\
\frac{1}{2}(W_i+Y_i)\mapsto E_{n+1,i}-E_{i,n+1}\qquad\text{for}\quad 1\leq i\leq n,
\\
\frac{1}{2}(W_i-Y_i)\mapsto E_{n+2,i}+E_{i,n+2}\qquad\text{for}\quad 1\leq i\leq n,
\\
Z_{i,j}\mapsto E_{i,j}-E_{j,i}\qquad\text{for}\quad 1\leq i<j\leq n,
\end{gather*}
where $E_{k,l}$ is the $(n+2)\times(n+2)$ matrix with single non-zero entry in
the $k$th row and~$l$th column.
The generators $X_1$ and~$X_3$ generate the parabolic subalgebra of $\mathfrak{sl}_2(\R)$ of upper
triangular matrices.

\section{The group}
\label{Three}

With an eye toward the construction of a~family of induced representations of the group
$G:=\text{SL}(2,\R)\times\text{SO}(n+1,1)_0$, we will study its structure in this section.
Since most of this material is standard, most details will be omitted
(for a~general treatment see~\cite[Chapter~VI]{Knapp01}).

Let $\mathfrak{g}=\mathfrak{k}\oplus\mathfrak{p}$ be the Cartan decomposition under
the standard Cartan involution.
Then, $\mathfrak{k}$ is isomorphic to $\mathfrak{so}(2)\times\mathfrak{so}(n+1)$
with $\mathfrak{so}(2)\subset\mathfrak{sl}_2(\R)$ and
$\mathfrak{so}(n+1)$ embedded in $\mathfrak{so}(n+1,1)$ in the upper left $(n+1)\times(n+1)$ block.

We consider the minimal parabolic subalgebras $\mathfrak q^\pm$ of $\mathfrak g$ and~their
respective Langlands decompositions $\mathfrak{m\oplus a\oplus n^\pm}$.
The maximal Abelian subalgebra $\mathfrak{a}\subset\mathfrak{p}$ is given by
\begin{gather*}
\mathfrak{a}:=\text{span}\left\{H_{v,y}:=\left.\left[\begin{pmatrix}
v&0
\\
0&-v
\end{pmatrix},
\left(\begin{array}{@{}c|c@{}}
0_{(n+1)\times(n+1)}&\begin{array}{@{}c@{}}0
\\
\vdots
\\
0
\\
-y\end{array}
\\
\hline
\begin{matrix}0&\cdots&0&-y\end{matrix}&0
\end{array}\right)\right]\;\right|\;v,y\in\R\right\}.
\end{gather*}
The centralizer of $\mathfrak{a}$ in $\mathfrak{k}$, that is denoted by $\mathfrak m$, is trivial
in the $\mathfrak{sl}_2(\R)$ component and~it embeds as $\mathfrak{so}(n)$ in the upper left
$n\times n$
block in the $\mathfrak{so}(n+1,1)$ component.

If
$\nu^\pm_{i}:=E_{n+1,i}-E_{i,n+1}\pm E_{n+2,i}\pm E_{i,n+2}\in\text{Mat}_{(n+2)\times(n+2)}$
for $1\leq i\leq n$, then the nilpotent subalgebras
\begin{gather*}
\mathfrak n^+=\spn\left\{\eta^+_{i,a,\sigma}:=\left[\begin{pmatrix}
0&\sigma
\\
0&0
\end{pmatrix},a\nu^+_i\right]\;\Bigg|\;a,\sigma\in\R,1\leq i\leq n\right\}
\end{gather*}
and
\begin{gather*}
\mathfrak n^-=\spn\left\{\eta^-_{i,a,\sigma}:=\left[\begin{pmatrix}
0&0
\\
\sigma&0
\end{pmatrix},a\nu^-_i\right]\;\Bigg|\;a,\sigma\in\R,1\leq i\leq n\right\}.
\end{gather*}
Consequently, the
minimal parabolic subalgebras of $\mathfrak g$ are def\/ined as $\mathfrak q=\mathfrak{m\oplus
a\oplus n}$ and
$\mathfrak q^-=\mathfrak{m\oplus a\oplus n^-}$.

At the level of the group, these subalgebras exponentiate to
\begin{gather*}
A:=\left\{h_{a,y}:=\left.\left[\begin{pmatrix}
a&0
\\
0&a^{-1}
\end{pmatrix},
\left(\begin{array}{@{}c|c@{}}
I_{n}&0
\\
\hline
0&\begin{matrix}
\cosh y&\sinh y
\\
\sinh y&\cosh y
\end{matrix}
\end{array}\right)\right]\;\right|\;y\in\R,a>0\right\},
\\
M:=\left\{m_{j,B}:=\left[\begin{pmatrix}
-1&0
\\
0&-1
\end{pmatrix}^j,
\left(\begin{array}{@{}c|c@{}}
B&0
\\
\hline
0&I_2
\end{array}\right)\right]\;\Bigg|\;B\in {\rm SO}(n),j\in\mathbb Z_2\right\},
\\
N:=\left\{n_{t,x}:=\left.\left[\begin{pmatrix}
1&t
\\
0&1
\end{pmatrix},
\left(\begin{array}{@{}c|c@{}}
I_{n}&\begin{matrix}
-x&&&x
\end{matrix}
\\
\hline
\begin{matrix}
x^\text{T}
\\
x^\text{T}
\end{matrix}
&\begin{matrix}
1-\frac{1}{2}\|x\|^2&\frac{1}{2}\|x\|^2\tsep{3pt}
\\
-\frac{1}{2}\|x\|^2&1+\frac{1}{2}\|x\|^2\tsep{2pt}
\end{matrix}
\end{array}\right)\right]\;\right|\;(t,x)\in\R^{1,n}\right\},
\end{gather*}
and
\begin{gather*}
N^-:=\left\{n_{t,x}^-:=\left.\left[\begin{pmatrix}
1&0
\\
t&1
\end{pmatrix},
\left(\begin{array}{@{}c|c@{}}
I_{n}&\begin{matrix}
x&&&&x
\end{matrix}
\\
\hline
\begin{matrix}
-x^\text{T}
\\
x^\text{T}
\end{matrix}
&\begin{matrix}1-\frac{1}{2}\|x\|^2&-\frac{1}{2}\|x\|^2\tsep{3pt}
\\
\frac{1}{2}\|x\|^2&1+\frac{1}{2}\|x\|^2\tsep{2pt}
\end{matrix}
\end{array}\right)\right]\;\right|\;(t,x)\in\R^{1,n}\right\}.
\end{gather*}
The minimal parabolic subgroups corresponding to $\mathfrak{q}$ and~$\mathfrak{q}^-$ are $Q=MAN$
and $Q^-=MAN^-$ respectively.

\section{Induced representations}
\label{induced}
In this section, we will construct a~family of characters on $Q^-$.
We will induce a~representation of $G$ from each of these characters and~we will identify each of
these induced representations with a~subspace of $C^\infty(\R^{1,n})$.
It is shown in~\cite{Vogan} that $G$ acts globally on these special classes of functions.
Since the symmetry group of the porous medium equation can be embedded in $G$, the action of $G$
restricts to a~global action of the symmetry group of the porous medium equation.

We start by noticing that a~general element $\phi_{b,c}\in\mathfrak{a}^*$ is a~linear functional
that acts by
$\phi_{b,c}(H_{v,y})=bv+cy$ for f\/ixed constants $b,c\in\C$.
By exponentiating $\phi_{b,c}$ we obtain
a character on $A$.
The resulting character, $\chi_{r,s}:A\to\C$, is determined by two continuous parameters $r,s\in\C$
and def\/ined by
\begin{gather*}
\chi_{r,s}(h_{a,y})=a^r e^{sy}.
\end{gather*}
A character $\chi_p:M\to\C$ can be def\/ined by
\begin{gather*}
\chi_{p}(m_{j,B})=(-1)^{jp},
\end{gather*}
for $p\in\mathbb Z_2$.
We will consider the characters parametrized by two continuous parameters $r,s\in\C$ and
a discrete parameter $p\in\mathbb Z_2$, $\chi_{p,r,s}:Q^-\to\C$ def\/ined by
\begin{gather*}
\chi_{p,r,s}(q^-)=\chi_p(m)\chi_{r,s}(a),
\end{gather*}
where $q^-=man^-$ and~by requiring it to be trivial on $\mathfrak{n}^-$.

We consider the inf\/inite dimensional induced representation
\begin{gather*}
\text{Ind}_{Q^-}^G(\chi_{p,r,s}):=\left\{\varphi\in C^\infty(G)\;\Big|\;\varphi(gq^-)
=\chi_{p,r,s}(q^-)^{-1}\varphi(g)\;\text{for}\; g\in G\;\text{and}\;q^-\in Q^-\right\},
\end{gather*}
with the $G$-action def\/ined by left translation.
This space is known as the induced picture.

The unipotent group $N$ is isomorphic to $\R^{1,n}$ via $(t,x)\mapsto n_{t,x}$.
Since $N Q^-$ embeds in $G$ as an open dense set, it is easy to see from the def\/inition of
$\text{Ind}_{Q^-}^G(\chi_{p,r,s})$ that an element $\varphi\in\text{Ind}_{Q^-}^G(\chi_{p,r,s})$ is
completely determined by its value on $N$.
Therefore, the restriction map between $\text{Ind}_{Q^-}^G(\chi_{p,r,s})$ and~the space
\begin{gather}
\label{Iprime}
I'(p,r,s):=\left\{f\in C^\infty(\R^{1,n})\;\Big|\;f(t,x)=\varphi(n_{t,x})
\;\text{for some}\;\varphi\in\text{Ind}_{Q^-}^G(\chi_{p,r,s})\right\}
\end{gather}
is injective.
By the def\/inition of $I'(p,r,s)$, the restriction map $\varphi\to f$ is surjective, thus an
isomorphism of vector spaces.
A $G$-module structure can be given to $I'(p,r,s)$ so that the map $\varphi\to f$ is intertwining.
Therefore, $\text{Ind}_{Q^-}^G(\chi_{p,r,s})$ is isomorphic to $I'(p,r,s)$ as $G$-module.
This space is known as the non-compact picture.

\section[Actions on $I'(p,r,s)$]{Actions on $\boldsymbol{I'(p,r,s)}$}
\label{actions}

In this section, we will describe the actions of $G$ and~of $\mathfrak{g}$ on $I'(p,r,s)$.
As a~result, we will f\/ind special values of $r$ and~$s$ that will determine the class of
functions on which the action of the symmetry group of equation~\eqref{Porous}, globalizes.

\subsection[Actions of $\mathfrak{sl}(2,\R)$ and~$\text{SL}(2,\R)$]{Actions of $\boldsymbol{\mathfrak{sl}(2,\R)}$ and~$\boldsymbol{\text{SL}(2,\R)}$}

In this section, when we consider the action of an element $g\in\text{SL}(2,\R)$, it should be
understood as the action of its image under the natural inclusion map
$\text{SL}(2,\R)\hookrightarrow G$.
The analogue will be true for elements in the Lie algebra.
If possible, we will avoid writing this image explicitly to avoid cumbersome notation.
\begin{proposition}
\label{uno}
An element $\left(\begin{smallmatrix}a~&b
\\
0&a^{-1}\end{smallmatrix}\right)\in\text{\rm SL}(2,\R)$ acts on $I'(p,r,s)$ by
\begin{gather*}
\left(\begin{matrix}a&b
\\
0&a^{-1}\end{matrix}\right).f(t,x)=\sgn(a)^{p}|a|^r f\left(\frac{t-ab}{a^2},x\right).
\end{gather*}
Let $\{H,E,F\}$ be the standard basis for $\mathfrak{sl}(2,\R)$.
Then, $H$ acts by the differential operator $-2t\partial_t+r u\partial_u$
and $E$ acts by $\partial_t$
on $I'(p,r,s)$.
\end{proposition}

\begin{proof}
First, we write
\[
\left(\begin{matrix}a &b
\\
0&a^{-1}\end{matrix}\right)=\left(\begin{matrix}1&b a
\\
0&1\end{matrix}\right)\left(\begin{matrix}a &0
\\
0&a^{-1}\end{matrix}\right)
\]
 and~we calculate each action of each element individually.
Let $f\in I'(p,r,s)$.
Then, using the def\/inition of $I'(p,r,s)$ and~the $G$ action on it, we have
\begin{gather*}
\left(\begin{matrix}a&0
\\
0&a^{-1}\end{matrix}\right).f(t,x)=\left(\begin{matrix}a&0
\\
0&a^{-1}\end{matrix}\right).\varphi(n_{t,x})=\varphi\left(\left[\left(\begin{matrix}a^{-1}&0
\\
0&a\end{matrix}\right),I_{n+2}\right]n_{t,x}\right).
\end{gather*}
Now we notice that
\begin{gather*}
\left[\left(\begin{matrix}a^{-1}&0
\\
0&a\end{matrix}\right),I_{n+2}\right]n_{t,x}=
n_{a^{-2}t,x}\cdot\left[\left(\begin{matrix}a^{-1}&0
\\
0&a\end{matrix}\right),I_{n+2}\right].
\end{gather*}
Using the character $\chi_{p,r,s}$ and~the def\/inition of $\text{Ind}_{Q^-}^G(\chi_{p,r,s})$, we
obtain
\begin{gather*}
\left(\begin{matrix} a & 0
\\
0 & a^{-1}
\end{matrix}\right).f(t,x)=\chi_{p,r,s}\left(\left[(\sgn(a)\left(\begin{matrix} |a|^{-1}
& 0
\\
0 & |a| \end{matrix}\right),I_{n+2}\right]\right)^{-1} \varphi( n_{a^{-2}t, x})\\
\hphantom{\left(\begin{matrix} a & 0
\\
0 & a^{-1}
\end{matrix}\right).f(t,x)}{}
=\sgn(a)^p|a|^r f( a^{-2}t, x).
\end{gather*}
The calculation for the other element is similar, thus omitted.
The result on the Lie algebra action follows by dif\/ferentiation.
\end{proof}

\begin{remark}
Setting $r=\frac{2}{m-1}$ we recover the actions of the elements $X_1$ and~$X_3$ in the symmetry
group.
We did not calculate the action of $F$ because the space of solutions to the porous medium equation
is not invariant under the subalgebra generated by this element.
So, for the goals of this paper, the action of this element is not relevant.
\end{remark}
\begin{remark}
Notice that exponentiation of the inf\/initesimal generators in $\mathfrak{sl}(2,\R)$ gives the
action of $\left(\begin{smallmatrix}a~&b
\\
0&a^{-1}\end{smallmatrix}\right)$ when $a>0$.
With the appropriate value of $r$ this agrees with the action calculated in the previous
proposition.
However, the previous proposition gives us a~way to extend the action when $a<0$.
This ref\/lects the global nature of the action of $G$ on $I'(p,r,s)$.
\end{remark}

\subsection[Actions of $\mathfrak{so}(n+1,1)$ and~$\text{SO}(n+1,1)_0$]{Actions of $\boldsymbol{\mathfrak{so}(n+1,1)}$ and~$\boldsymbol{\text{SO}(n+1,1)_0}$}

As in the previous section, in this section when we consider the action of an element
$g\in\text{SO}(n+1,1)_0$, it should be understood as the action of its image under the natural
inclusion map $\text{SO}(n+1,1)_0\hookrightarrow G$.
The same convention will be used for the actions of the Lie algebra.

\subsubsection[Actions of $\mathfrak{m}$ and~$M$]{Actions of $\boldsymbol{\mathfrak{m}}$ and~$\boldsymbol{M}$}
\begin{lemma}
Let $g_{i,j,\theta}=\exp_{{\rm SO}(n+1,1)}(\theta(E_{i,j}-E_{j,i}))$ for $1\leq i<j\leq n$. Then
\begin{gather*}
[I_2,g_{i,j,-\theta}]\cdot n_{t,x}=n_{t,g_{i,j,-\theta}x}\cdot[I_2,g_{i,j,-\theta}].
\end{gather*}
\end{lemma}

\begin{proof}
This is a~straightforward matrix calculation.
\end{proof}

\begin{proposition}
An element $g_{i,j,\theta}$ for $1\leq i<j\leq n$ acts on a~function $f\in I'(p,r,s)$ by
\begin{gather*}
g_{i,j,\theta}.f(t,x)=f(t,g_{i,j,-\theta}x).
\end{gather*}
The element $E_{i,j}-E_{j,i}\in\mathfrak{so}(n+1,1)$ for $1\leq i<j\leq n$ acts on $I'(p,r,s)$ by
the differential operator
$
x_i\partial_j-x_j\partial_i$.
\end{proposition}

\begin{proof}
By the def\/inition of $I'(p,r,s)$ and~by the action of $G$ on $\text{Ind}_{Q^-}^G(\chi_{p,r,s})$
we have
\begin{gather*}
g_{i,j,\theta}.f(t,x)=g_{i,j,\theta}.\varphi(n_{t,x})=\varphi([I_2,g_{i,j,-\theta}].
n_{t,x}).
\end{gather*}
The result now follows from the lemma and~the def\/inition of $\chi_{p,r,s}$.
The result for the Lie algebra follows by taking $\frac{d}{d\theta}|_{\theta=0}$.
\end{proof}

\subsubsection[Actions of $\mathfrak{n}$ and~$N$]{Actions of $\boldsymbol{\mathfrak{n}}$ and~$\boldsymbol{N}$}

A straightforward matrix calculation shows that
\begin{gather*}
n_{t',x_i e_i}n_{t,x}=n_{t-t',x-x_i e_i},
\end{gather*}
where $e_i$ is the standard basis element of $\R^n$.
We use this fact to prove the following proposition.
\begin{proposition}
The action of $n_{t',x'}$ on $f\in I'(p,r,s)$ is given by
\begin{gather*}
n_{t',x'}.f(t,x)=f(t-t',x-x')
\end{gather*}
and $\nu_i^+$ acts on $I'(p,r,s)$ by the differential operator
$\partial_i$, for $1\leq i\leq n$.
\end{proposition}

\subsubsection[Actions of $\mathfrak{a}$ and~$A$]{Actions of $\boldsymbol{\mathfrak{a}}$ and~$\boldsymbol{A}$}

Using the fact that
$
h_{1,\epsilon}^{-1}n_{t,x}=n_{t,e^{\epsilon}x}h_{1,\epsilon}^{-1}$,
we prove the following proposition.
\begin{proposition}
The action of $h_{1,\epsilon}$ on $f\in I'(p,r,s)$ is given by
$
h_{1,\epsilon}.f(t,x)=e^{s\epsilon}f\left(t,e^{\epsilon}x\right)
$
and~$H_{0,1}$ acts on~$I'(p,r,s)$ by the differential operator
$\sum\limits_{i=1}^n x_i\partial_i+su\partial_u$.
\end{proposition}

\begin{remark}
When $s=\frac{2}{m-1}$, this action corresponds to the action of the inf\/initesimal genera\-tor~$X_2$ of the symmetry group of the porous medium equation~\eqref{Porous}.
\end{remark}

\subsubsection[Actions of $\mathfrak{n}^-$ and~$N^-$]{Actions of $\boldsymbol{\mathfrak{n}^-}$ and~$\boldsymbol{N^-}$}

To describe this action, we need to introduce the maps $\delta_i:\R^n\to\R$ def\/ined by
\begin{gather*}
\delta_i(x)=1-2x_i+\|x\|^2.
\end{gather*}
We will also need the maps $\gamma_i:\R^{1,n}\to\R^n$ def\/ined by
\begin{gather*}
\gamma_i(\epsilon,x)=\delta_i(\epsilon x)^{-1}\big(x-\epsilon\|x\|^2e_i\big)
\end{gather*}
and the maps $\kappa_i:\R^{1,n}\to\R^n$ given by
\begin{gather*}
\kappa_i(\epsilon,x)=\delta_i(\epsilon x)^{-1}\epsilon(\epsilon x-e_i).
\end{gather*}
Then
\begin{gather*}
\exp\big(\epsilon\nu_i^-\big).f(t,x)=n^{-}_{0,\epsilon e_i}.\varphi\big(n_{t,x}\big)=
\varphi\big(n^{-}_{0,-\epsilon e_i}.n_{t,x}\big).
\end{gather*}
Now, in order to write the action back in terms of $f$ we need to decompose $n^{-}_{0,-\epsilon
e_i}.n_{t,x}$ as a~product of its $N\times M A N^-$ components.
\begin{lemma} For some $m\in M$
\begin{gather*}
n^{-}_{0,-\epsilon e_i}.n_{t,x}=
n_{t,\gamma_i(\epsilon,x)}m h_{1,-\log(\delta_i(\epsilon x))}n^-_{0,\kappa_i(\epsilon,x)}.
\end{gather*}
\end{lemma}

\begin{proof}
Since $\text{SO}(n+1,1)$ has real rank $1$, the Weyl group generated by the restricted roots
$W(G,A)$ has two elements.
The non-trivial element $\omega\in W(G,A)$ acts on $N$ by $\omega^{-1}N\omega=N^-$.
Then, by the uniqueness of the Bruhat decomposition, it suf\/f\/ices to show that
\begin{gather*}
\big(n_{t,\gamma_i(\epsilon,x)}\big)^{-1}n^{-}_{0,-\epsilon e_i}n_{t,x}
\big(n^-_{0,\kappa_i(\epsilon,x)}\big)^{-1}h_{1,-\log(\delta_i(\epsilon x))}^{-1}
\end{gather*}
dif\/fers from the $(n+2)\times(n+2)$ identity matrix only in the upper left $n\times n$ block
(i.e.\ it is in~$M$).
This is a~straightforward but long matrix calculation.
\end{proof}

\begin{proposition}
\label{cinco}
Let $f\in I'(p,r,s)$.
Then,
\begin{gather*}
n^{-}_{0,-\epsilon e_i}f(t,x)=\delta_i(\epsilon x)^{-s}f(t,\gamma_i(\epsilon,x))
\end{gather*}
and $\nu_i^-$ acts on $I'(p,r,s)$ as the differential operator
\begin{gather*}
\left(x_i^2-\sum_{j\neq i}x_j^2\right)\partial_i+\sum_{j\neq i}2x_ix_j\partial_j+2x_is u\partial_u.
\end{gather*}
\end{proposition}

\begin{proof}
The f\/irst equation follows from the previous lemma and~the second equation follows form
dif\/ferentiating the f\/irst.
\end{proof}

We summarize the previous propositions in the following theorem.
\begin{theorem}
The action of the group $G$ on $I'\left(p,\frac{2}{m-1},\frac{2}{m-1}\right)$, for $p\in\mathbb
Z_2$, gives a~globalization of the action of the local symmetry group of equation~\eqref{Porous}.
\end{theorem}

\begin{proof}
Follows from Propositions~\ref{uno} to~\ref{cinco}.
\end{proof}

\section{Applications}
\label{applications}

The problem of explicitly describing the solution space of the porous medium equation inside
$I'\left(p,\frac{2}{m-1},\frac{2}{m-1}\right)$ is a~complicated problem.
A possible strategy is to look at the compact picture of
$I'\left(p,\frac{2}{m-1},\frac{2}{m-1}\right)$ and~analyze what conditions need to be satisf\/ied
by the functions in this space to correspond to a~solution of the dif\/ferential equation.
For references in this direction see~\cite{Franco} and~\cite{Sepanski3}.
However, some solutions can be realized in $I'\left(p,\frac{2}{m-1},\frac{2}{m-1}\right)$ and~for
illustration purposes we will examine a~few in this section.

\subsection{Stationary solutions}
Since we work in a~smooth class, the simplest example to consider is provided by stationary
solutions.
In this case, the solutions are given as harmonic polynomials raised to the $1/m$th power.

Let $k:\R^n\to\R$ be a~harmonic polynomial and~def\/ine $f\in C^\infty(\R^{1,n})$ by
$f(t,x)=k(x)^{1/m}$.
It is well-known that $f$ satisf\/ies equation~\eqref{Porous}.

Since sections in $I\left(p,\frac{2}{m-1},\frac{2}{m-1}\right)$ are completely determined by their
values on $N$ and~we know the values that $\varphi$ must take, namely $\varphi(n_{t,x})=f(t,x$), we
can extend $\varphi$ to $N MAN^-$ using the character $\chi_{p,\frac{2}{m-1},\frac{2}{m-1}}$.
This can be used to def\/ine $\varphi$ on all of $G$ via limits, because $N MAN^-$ sits as an open
dense subset of $G$.
To exemplify how this process works, let us choose the simplest example.
Let $k:\R\to\R$ be given by $k(x)=x$.
By explicitly decomposing the element
\begin{gather*}
\left[
\begin{pmatrix}a&b
\\
c&d\end{pmatrix},\begin{pmatrix}a_{11}&a_{12}&a_{13}
\\
a_{21}&a_{22}&a_{23}
\\
a_{31}&a_{32}&a_{33}\end{pmatrix}\right]
\end{gather*}
in $G$ in its $N\times MAN^-$ components, we determine
that
\begin{gather*}
\varphi
\left(
\left[
\begin{pmatrix} a& b
\\
c& d \end{pmatrix},\begin{pmatrix}a_{11} & a_{12} & a_{13}
\\
a_{21} & a_{22} & a_{23}
\\
a_{31} & a_{32} & a_{33}
\end{pmatrix}
\right]\right)
\\
\qquad\qquad
=\sgn(d)^p|d|^{\frac{2}{m-1}}  \left(\frac{a_{21}+a_{31}}{1+a_{11}}\right)^{\frac{1}{m}}
\left(\frac{-2(a_{21}+a_{31})}{(1+a_{11})(a_{12}-a_{13})}\right)^{\frac{2}{1-m}}
\end{gather*}
for $a_{11}\neq-1$ and~$a_{12}\neq a_{13}$.
As expected, $\varphi$ restricts to $f$ on $N$, more specif\/ically $\varphi(n_{t,x})=f(t,x)$.
To extend $\varphi$ to the elements in $G$ for which $a_{11}=-1$, we notice that in this case the
conditions on the group elements force $a_{12}=\pm a_{13}$ and~$a_{21}=\pm a_{31}$.
If $a_{11}=-1$, $a_{12}=\pm a_{13}\neq0$, and~$a_{21}=-a_{31}$, we can use limits to determine the
appropriate value for $\varphi$
\begin{gather*}
\varphi\left(\left[\begin{pmatrix}a&b
\\
c&d\end{pmatrix},\begin{pmatrix}-1&\pm a_{13}&a_{13}
\\
-a_{31}&a_{22}&a_{23}
\\
a_{31}&a_{32}&a_{33}\end{pmatrix}\right]\right)=
\sgn(d)^p|d|^{\frac{2}{m-1}}
\left(\frac{\mp1}{a_{31}}\right)^{\frac{1}{m}}\left(\frac{\mp2}{a_{31}a_{13}}\right)^{\frac{2}{1-m}}.
\end{gather*}
The values of $\varphi$ when $a_{21}=a_{31}$ and~when $a_{12}=a_{13}=0$ can be determined in
a~similar fashion.
The resulting map can be shown to be smooth and~it was constructed in such way that the condition
$\varphi(gq^-)=\chi_{p,\frac{2}{m-1},\frac{2}{m-1}}(q^-)^{-1}\varphi(g)$ is satisf\/ied.
Therefore, the so def\/ined $\varphi$ belongs to $I\left(p,\frac{2}{m-1},\frac{2}{m-1}\right)$
and~its image $f$ belongs to $I'\left(p,\frac{2}{m-1},\frac{2}{m-1}\right)$.

This procedure can be repeated for higher dimensions.
For example, if we let $k:\R^2\to\R$ be a~2-dimensional harmonic polynomial, then the corresponding
map $\varphi$ is given by
\begin{gather*}\varphi\left(\left[\begin{pmatrix} a& b
\\
c& d \end{pmatrix},\begin{pmatrix}a_{11} & a_{12} & a_{13} & a_{14}
\\
a_{21} & a_{22} & a_{23} & a_{24}
\\
a_{31} & a_{32} & a_{33} & a_{34}
\\
a_{41} & a_{42} & a_{43} & a_{44} \end{pmatrix}\right]\right)
=\sgn(d)^p|d|^{\frac{2}{m-1}}(k(z_1,z_2))^{\frac{1}{m}}
\\
\qquad\qquad
\times
\big(a_{44}\big(1+z_1^2+z_2^2\big)-a_{34}\big(-1+z_1^2+z_2^2\big)-2(a_{14}z_1+a_{24}z_2)\big)^\frac{2}{1-m},
\end{gather*}
where
\begin{gather*}
z_1=\frac{(a_{11}+a_{22})(a_{31}+a_{41})-(a_{12}-a_{21})(a_{32}+a_{42})}{(a_{12}-a_{21})^2+(a_{11}+a_{22})^2}
\end{gather*}
and
\begin{gather*}
z_2=\frac{(a_{11}+a_{22})(a_{32}+a_{42})-(a_{12}-a_{21})(a_{31}+a_{41})}{(a_{12}-a_{21})^2+(a_{11}+a_{22})^2}
\end{gather*}
whenever $a_{12}\neq a_{21}$ or $a_{11}\neq-a_{22}$.
In this case determining the smooth extension to the whole group is more involved, but it mimics
the procedure used in the one-dimensional case.
It is easy to see that the so def\/ined map $\varphi$ restricts to $k(x)^{1/m}$ on $N$ as desired.

\subsection{Other solutions}
To stay within the smooth class of functions, some of the well-known solutions of the porous medium
equation can only be considered for specif\/ic values of $m$.
This is the case of the solutions obtained via separation of variables.
It is well known that these solutions are of the form
\begin{gather*}
u(t,x)=((m-1)(t-t_0))^{-1/(m-1)}F(x),
\end{gather*}
where $\Delta F^m(x)+F(x)=0$, see~\cite[Chapter~4]{vazquez2006porous}.
When $m>1$ these solutions have a~non-re\-mo\-vable singularity at $t=t_0$.
However, for $m<1$ these solutions satisfy the condition $u(t_0,x)=0$.
Solutions of this type would need to be extended to sections in
$I\big(p,\frac{2}{m-1},\frac{2}{m-1}\big)$ in a~similar fashion as we extended the stationary
solutions.
This in itself can be a~complicated problem.

For source-type solutions, the smooth category would need to be abandoned and~other categories may
be considered.
However, the induced representation may no longer carry the structure of a~global $G$-module.
For more information in this direction see~\cite{Vogan}.

\section{Compact picture}

To study the structure of $\text{Ind}_{Q^-}^G\Big(\chi_{p,\frac{2}{m-1},\frac{2}{m-1}}\Big)$ as
a~$G$-representation, it is useful to look at an isomorphic copy of it called the compact picture.
To construct $I'\Big(p,\frac{2}{m-1},\frac{2}{m-1}\Big)$, we used restriction to the non-compact
subgroup $N\cong\R^{1,n}$.
To construct the compact picture of $\text{Ind}_{Q^-}^G\Big(\chi_{p,\frac{2}{m-1},\frac{2}{m-1}}\Big)$ we
use restriction to the maximal compact subgroup $K\subset G$.

The Iwasawa decomposition of $G$ is given by $G=KAN^-$.
Since $AN^-\subset Q^-$, a~map in the induced picture is completely determined by its restriction
to $K$.
The compact picture is def\/ined as the image of this restriction, that is
\begin{gather*}
I''\left(p,\frac{2}{m-1},\frac{2}{m-1}\right)=\Big\{ \zeta \in C^\infty(K)
\,\Big| \, \exists \, \varphi \in \text{Ind}_{Q^-}^G\Big(\chi_{p,\frac{2}{m-1},\frac{2}{m-1}}\Big):
\zeta(k)=\varphi(k)\; 
\forall\, k\in K \Big\}
\end{gather*}
and it is isomorphic to $\text{Ind}_{Q^-}^G\Big(\chi_{p,\frac{2}{m-1},\frac{2}{m-1}}\Big)$ as vector spaces.
The space can be given the structure of a~$G$-module so that the restriction map is intertwining.
Consequently,
\begin{gather*}
\text{Ind}_{Q^-}^G\Big(\chi_{p,\frac{2}{m-1},\frac{2}{m-1}}\Big)\cong I'\left(p,\frac{2}{m-1},\frac{2}{m-1}\right)
\end{gather*}
and~$I'\left(p,\frac{2}{m-1},\frac{2}{m-1}\right)\cong
I''\left(p,\frac{2}{m-1},\frac{2}{m-1}\right)$ as $G$-modules.
The space
$I''\left(p,\frac{2}{m-1},\frac{2}{m-1}\right)$ contains an open dense subset given by
\begin{gather*}
\Big\{\zeta\in C^\infty(K)\,\Big|\,\zeta(k\mu)=
\chi_{p,\frac{2}{m-1},\frac{2}{m-1}}(\mu)^{-1}\zeta(k)\;\text{for}\; k\in K\;\text{and}\;\mu\in M
\cap K\Big\}
\end{gather*}
(see~\cite[Chapter~2]{knapp2001representation}).
Since, $\text{SO}(n+1)/\text{SO}(n)\cong S^{n}$ we obtain
\begin{gather*}
\Big\{F\in C^\infty(S^1\times S^{n})\,\Big|\,F(\theta+j\pi,z)=
(-1)^{jp}F(\theta,z)\;\text{for}\; j\in\mathbb Z\Big\}
\end{gather*}
is an open dense set in $I''\left(p,\frac{2}{m-1},\frac{2}{m-1}\right)$.
This gives a~more concrete realization of the representation spaces and~exhibits their inf\/inite
dimensionality.

\pdfbookmark[1]{References}{ref}
\LastPageEnding

\end{document}